\documentclass[10pt]{amsart}
\usepackage{graphicx}
\usepackage{amssymb,amscd,xcolor}
\usepackage{tikz}
\usetikzlibrary{cd}
\usepackage{caption}

\usepackage{ifpdf}
\ifpdf\usepackage[colorlinks,final,hyperindex]{hyperref}
\else\usepackage[colorlinks,final,hyperindex,hypertex]{hyperref}
\fi
\usepackage[backend=biber,style=alphabetic,hyperref,backref,sorting=nyt]{biblatex}
\newtheorem{thm}{Theorem}[section]
\newtheorem{defi}[thm]{Definition}
\newtheorem{rmk}[thm]{Remark}
\newtheorem{lem}[thm]{Lemma}
\newtheorem{prop}[thm]{Proposition}
\newtheorem{cor}[thm]{Corollary}
\newtheorem{exa}[thm]{Example}

\title{Functoriality of toric coherent-constructible correspondence}
\author{Yuze Sun}
\address{Beijing International Center for Mathematical Research, Peking University, 5 Yiheyuan Road, Beijing 100871, China}
\email{sunyuze@stu.pku.edu.cn}
\date{\today}

\begin{document}

\newcommand{\Vect}{\mathsf{Vect}}
\newcommand{\Spec}{\operatorname{Spec}}
\newcommand {\QCoh}{\mathsf{QCoh}}
\newcommand{\indcoh}{\mathsf{IndCoh}}
\newcommand {\Coh}{\mathsf{Coh}}
\newcommand {\Perf}{\mathsf{Perf}}
\newcommand {\Sh}{\mathsf{Sh}}
\newcommand {\Fun}{\mathsf{Fun}}
\newcommand {\op}{\operatorname{op}}
\newcommand {\Mod}{\mathsf{Mod}}
\newcommand{\coMod}{\mathsf{coMod}}
\renewcommand{\ss}{\operatorname{ss}}
\newcommand {\prlst}{\mathsf{Pr}^{\L}_{\mathsf{st}}}
\newcommand {\prlstk}{\mathsf{Pr}^{\L}_{\mathsf{st},k}}
\newcommand{\ex}{\operatorname{ex}}
\newcommand{\Hom}{\operatorname{Hom}}
\renewcommand{\hom}{\mathcal{H}om}
\newcommand{\catex}{\mathsf{Cat}^{\operatorname{Ex}}}
\renewcommand{\dim}{\operatorname{dim}}
\newcommand{\cech}{\v{C}ech}
\newcommand{\spec}{\operatorname{Spec}}
\newcommand{\calg}{\operatorname{CAlg}}
\newcommand{\rep}{\mathsf{Rep}}
\newcommand{\Int}{\operatorname{Int}}
\newcommand{\Star}{\operatorname{Star}}

\newcommand{\R}{\mathrm{R}}
\renewcommand{\L}{\mathrm{L}}
\newcommand{\lad}{\mathrm{LAd}}
\newcommand{\rad}{\mathrm{RAd}}
\newcommand {\lex}{\mathrm{LEx}}
\newcommand {\rex}{\mathrm{REx}}
\newcommand{\rank}{\mathrm{rank}}

\newcommand{\cA}{\mathcal{A}}
\newcommand{\cC}{\mathcal{C}}
\newcommand{\cD}{\mathcal{D}}
\newcommand{\cE}{\mathcal{E}}
\newcommand{\cF}{\mathcal{F}}
\newcommand{\cK}{\mathcal{K}}
\newcommand{\cM}{\mathcal{M}}
\newcommand{\cX}{\mathcal{X}}
\newcommand{\cO}{\mathcal{O}}

\newcommand{\bX}{{\mathbf{X}}}
\newcommand{\bU}{{\mathbf{U}}}
\newcommand{\bZ}{{\mathbf{Z}}}
\newcommand{\id}{\mathbf{id}}
\newcommand{\1}{\mathbf{1}}

\renewcommand{\lim}{\qopname\relax m{\mathbf{lim}}}
\newcommand{\colim}{\qopname\relax m{\mathbf{colim}}}
\newcommand{\fib}{\mathbf{fib}}
\newcommand{\cofib}{\mathbf{cofib}}
\renewcommand{\ker}{\mathbf{ker}}
\newcommand{\im}{\mathbf{im}}
\newcommand{\coker}{\mathbf{coker}}
\newcommand{\tot}{\mathbf{Tot}}
\newcommand{\Ab}{\mathbf{Ab}}
\newcommand{\bc}{\mathrm{Bar}}
\newcommand{\cb}{\mathrm{coBar}}
\newcommand{\oblv}{\mathbf{oblv}}

\newcommand{\GG}{\mathbb{G}}
\newcommand{\PP}{\mathbb{P}}
\newcommand{\RR}{\mathbb{R}}
\newcommand{\TT}{\mathbb{T}}
\newcommand{\ZZ}{\mathbb{Z}}
\renewcommand{\AA}{\mathbb{A}}
\newcommand{\QQ}{\mathbb{Q}}

\newcommand{\sN}{\mathsf{N}}
\newcommand{\sM}{\mathsf{M}}
\newcommand{\sL}{\mathsf{L}}
\newcommand{\sK}{\mathsf{K}}

\newcommand{\git}{\mathbin{
  \mathchoice{/\mkern-6mu/}
    {/\mkern-6mu/}
    {/\mkern-5mu/}
    {/\mkern-5mu/}}}

\begin{abstract}
    A morphism from a diagonalizable group $G$ to the torus of a toric variety $X$ induces an action of $G$ on $X$. We prove the category of ind-coherent sheaves on the quotient stack is equivalent to the category of sheaves on a cover of a real torus with singular supports contained in the FLTZ skeleton, extending Kuwagaki's nonequivariant coherent-constructible correspondence \cite{kuw16}. We also investigate the functoriality of such correspondence for toric morphisms and inclusions of orbit closures.
\end{abstract}

\maketitle

\section{Introduction}
As a variant of homological mirror symmetry, it was first observed by Bondal in \cite{Bondal} that the derived category of a toric variety $X_\Sigma$ is equivalent to a subcategory of constructible sheaves on $\sM_\RR/\sM$, where $\sM$ is the character lattice. Fang--Liu--Treumann--Zaslow \cite{FLTZ} proved $\mathsf{Perf}[X_\Sigma/T]\cong\Sh_{\Lambda}^{cc}(\sM_\RR)$ for any complete fan, where $\Sh^{cc}$ is the category of bounded complexes of compactly supported constructible sheaves. In \cite{FLTZ09}, a similar statement was proved for smooth complete toric DM stacks. In \cite{Tr10}, Treumann proved $\mathsf{Perf}(X_\Sigma)\cong\Sh^{cc}_\Lambda(\sM_\RR/\sM)$ for any smooth projective zonotopal fan, by showing the pullback map along $X\rightarrow[X_\Sigma/T]$ corresponds to the proper pushforward map along $\sM_\RR\rightarrow\sM_\RR/\sM$. This work was generalized to cragged stacky fans by Scherotzke--Sibilia in \cite{SS}. Kuwagaki \cite{kuw16} proved the nonequivariant equivalence $\indcoh(\cX_{\Sigma,\beta})\cong\Sh_\Lambda((\RR/\ZZ)^{\dim\Sigma})$ for toric stacks\footnote{The notion of toric stacks in \cite{kuw16} is from \cite{TorStk1} and is different from our notions in Sec.\ref{torstk}.}, see Thm.\ref{kuw}. Shende \cite{shende} proved $\indcoh[(\AA^n-Z)/G]\cong\Sh_\Lambda(\RR^n/\sM)$ for the Cox construction $(\AA^n-Z)\git G=X_{\Sigma}$ of a toric variety. Bai--Hu \cite{bai} proved $\QCoh[X_\Sigma/T]\cong\Sh_\Lambda(\sM_\RR)$ and $\QCoh(X_\Sigma)\cong\Sh_{\Lambda}(\sM_\RR/\sM)$ for any smooth projective fan over sphere spectrum coefficient. For the skeletons $\Lambda$'s above, see Def.\ref{skeleton}.

In this paper we extend the result of \cite{kuw16} to a larger class of toric stacks. Namely, we mean a quotient stack $[X_\Sigma/G]$ for a diagonalizable group $G$ acting on $X_{\Sigma}$ via a morphism $G\rightarrow T$, see Sec.\ref{torstk}. By applying (de-)equivariantization formalism, we construct a cover $T_\beta$ of $\sM_\RR/\sM$ and prove $$\indcoh(\cX_{\Sigma,\beta})\cong\Sh_{\Lambda_\beta}(T_\beta).$$
See Thm.\ref{main_thm_cat}.

This equivalence allows us to compare the six-functor formalism on both sides. In Thm.\ref{adj_of_inc} and Cor.\ref{exc_funs}, we express the comparison with left and right adjoints of the inclusion $\iota:\Sh_{\Lambda}(T_\beta)\rightarrow\Sh(T_\beta)$, which were studied in \cite{WrpShv}. This result leads to Thm.\ref{lftr}, which states that the exceptional pullback of toric stacks corresponds to $\1_{\Sigma,\beta}\star\phi_{T!}(-)$, where $\1_{\Sigma,\beta}$ is the convolutional unit of $\Sh_{\Lambda_\beta}(T_\beta)$ and $\phi_T$ is the induced map between $T_\beta$'s, generalizing \cite[Thm. 3.8]{FLTZ}. In Thm.\ref{r_fct}, we deduce that $(\phi^!)^\R$ corresponds to $\phi_T^!$ under a weakly semistable condition.

As a consequence, we investigate the functors associated to the inclusion of an orbit closure, see Sec.\ref{closedorbits}. Besides, the restriction functor $i^!$ was calculated as a different formula in \cite{shende-gammage} for a smooth fan. It remains unknown why these constructions are equivalent.

\subsection*{Acknowledgments}
I would like to thank my advisor Bohan Fang for inspirations and encouragements and Yixiao Li for many useful discussions. This work is partially supported by National Key R\&D Program of China 2023YFA1009803.

\section{Conventions}
    All categorical notions (categories, functors, algebra and module structures) refer to $(\infty,1)$-structures in the sense of \cite{HTT,HA}.

    Fix a commutative Noetherian regular ring $k$ over $\QQ$. Unless otherwise specified, notions in algebraic geometry (schemes, morphisms) are defined over $k$. We omit the notations $\L$ and $\R$ for derived functors. Instead, they refer to left and right adjoints of functors.
    
    For a scheme or stack $X$, the category $\QCoh(X)$ is the unbounded derived category of quasi-coherent sheaves. For the category $\indcoh(X)$ and relavent functors, we refer to \cite{vol1}.

    For the theory of toric varieties, we refer to \cite{Cox}. For a finitely generated abelian group $\sL$, we write $\sL^\vee:=\Hom(\sL,\ZZ)$ and $\sL_\RR=\sL\otimes_\ZZ\RR$. Suppose $\sN$ is a {\bf lattice} (finitely generated free abelian group), by a {\bf cone} on $\sN_\RR$ we mean a strongly convex rational polyhedral cone $\sigma\subseteq \sN_\RR$. For such a cone, we set $\RR\sigma$ to be the subspace of $\sN$ spanned by $\sigma$ and $\ZZ\sigma=\RR\sigma\cap\sN$. Set $\sM=\sN^\vee$, then $\sigma^\perp\subseteq\sM_\RR$ is the subspace of functions vanishing on $\sigma$. The dual of a cone is defined as $$\sigma^\vee:=\{\mathbf{m}\in\sM_\RR|\langle \mathbf{m},\mathbf{n}\rangle\geq0\mbox{ for all }\mathbf{n}\in\sigma\}.$$
    
    In the theory of constructible sheaves, we restrict to paracompact locally compact Hausdorff spaces with finite covering dimensions. For a space $X$ and a presentable stable category $\cC\in\prlst$, the category $\Sh(X,\cC)$ is the full subcategory of $\Fun(\operatorname{Open}_X^{\op},\cC)$ spanned by the functors $\cF$ satisfying the following conditions. For simplicity, set $\Sh(-)=\Sh(-,\Mod_k)$.
    \begin{enumerate}
        \item $\cF(\emptyset)=0$.
        \item The diagram is a pullback diagram for any open subsets $U,V$.\\
            \begin{tikzcd}
            \cF(U\cup V) \arrow[r] \arrow[d] & \cF(U) \arrow[d] \\
            \cF(V) \arrow[r]& \cF(U\cap V)\end{tikzcd}
        \item For a filtered poset $I$ and a functor $i\mapsto U_i, I\rightarrow \operatorname{Open}_X$, the natural map $\cF(\bigcup_I U_i)\rightarrow\lim_I \cF(U_i)$ is an equivalence.
    \end{enumerate}
    
    With some technical lemmas \cite{lemma,unbounded}, the results in \cite{KS90} are appropriately generalized to $(\infty,1)$-categories of unbounded complexes. Given a closed conic subset $\Lambda\subseteq T^*M$ of a cotangent bundle, the full subcategory of $\Sh(M)$ with the restriction $\ss(-)\subseteq\Lambda$ is denoted by $\Sh_\Lambda(M)$. By a {\bf skeleton} we mean a closed conic subanalytic subset of a cotangent bundle.

    For a finitely generated abelian group $\sM$, we define a group scheme $D(\sM)=\spec k[\sM]$. This is a fully faithful contravariant functor from finitely generated abelian groups to algebraic groups, whose essential image is called {\bf diagonalizable groups}. Restricted to this subcategory, the quasi-inverse is $\Hom_{\operatorname{Alg.gp}}(-,\GG_m)$. We denote both functors by $D(-)$.

\section{coherent-constructible correspondence}

\subsection{(De-)Equivariantization}

Here we briefly introduce a lemma from \cite{shende}, which was originally proved in \cite{1aff}. \footnote{Here we adopt the quasi-coherent sheaves to match the context of \cite{1aff}. In the later use, we replace them by ind-coherent sheaves at no risk, because the groups are smooth.}\\

We identify a discrete group $G$ with the group scheme $\bigsqcup_{G}\spec k$. In this subsection, we focus on two classes of group schemes:\begin{enumerate}
    \item Finitely generated abelian groups (disctete groups for short);
    \item Diagonalizable groups.
\end{enumerate}

\begin{defi}
    For a group $G$ as above with multiplication map $m$, the category $\QCoh(G)$ is a comonoidal category with comultiplication $m^*$, and we say $G$ acts on a category $\cC$ or $\cC$ is a $G$-module if $\cC\in\coMod_{\QCoh(G)}(\prlst)$. For categories $\cC$ and $\cC^\prime$ with $G$ actions, we write $$\Hom_G(\cC,\cC^\prime):=\lim_{n\in\Delta}[\Fun^{\L,\mathrm{Ex}}(\cC,\cC^\prime\otimes\QCoh(G)^{\otimes n})]\in\prlst$$
    for the $G$-equivariant left adjoint functors.
\end{defi}

Then we introduce the {\bf equivariantization} functor $(-)^G$ and {\bf de-equivariantization} functor $(-)_G$.

\begin{defi}  Let $G$ act on $\Vect:=\Mod_k$ trivially.\begin{enumerate}
    \item For a $G$-module $\cC$, set $\cC^G=\Hom_G(\Vect,\cC)$. Naturally $\cC^G$ is a module of $\Hom_G(\Vect, \Vect)\in\calg(\prlst)$. The equivalence $\Hom_G(\Vect,\Vect)=(\rep(G),\otimes)$ enhances $(-)^G$ to a functor $$\coMod_{\QCoh(G)}(\prlst)\rightarrow\Mod_{\rep(G)}(\prlst).$$
    \item Let $\rep(G)$ act on $\Vect$ via forgetful functor. For a $\rep(G)$-module $\cD$, define $\cD_G=\Vect\otimes_{\rep(G)}\cD$. The trivial action $\Vect\rightarrow\Vect\otimes\QCoh(G)$ induces a $G$-action on $\cD_G$ and enhances $(-)_G$ to a functor $$\Mod_{\rep(G)}(\prlst)\rightarrow\coMod_{\QCoh(G)}(\prlst).$$
\end{enumerate}
\end{defi}

In \cite[Thm. 10.4.4]{1aff}, Gaitsgory proved the following theorem, sometimes called ``categorical Pontryagin duality".

\begin{thm}
   The (de-)equivariantization functors $$(-)^G:\coMod_{\QCoh(G)}(\prlst)\leftrightarrows\Mod_{\rep(G)}(\prlst):(-)_G$$ are equivalences of categories and are mutually inverse.
\end{thm}

Then we introduce the Catier dual groups. The following lemma can be found in \cite[Ch.10-11]{1aff}.\footnote{This is the only reason we assume $k$ to be a $\QQ$-algebra. If this assumtion is dropped, we require every discrete group $G$ satisfies the condition that $|G_{tors}|$ is invertible in $k$.}
\begin{lem}\label{dualgrp}
    There is a symmetric monoidal equivalence $(\QCoh(D(G)),\star)\rightarrow(\rep(G),\otimes)$, which satisfies:\begin{enumerate}
        \item If $G$ is discrete, the convolution $\star$ is defined to be $m_*$;
        \item If $\varphi:G^\prime\rightarrow G$ is a morphism between discrete groups, then the forgetful map $\rep(G)\rightarrow\rep(G^\prime)$ corresponds to $D(\varphi)_*$ under the equivalence;
        \item If $G$ is diagonalizable, the convolution $\star$ is defined to be $(m^*)^\L$;
        \item If $\varphi:G^\prime\rightarrow G$ is a morphism between diagonalizable groups, then the forgetful map $\rep(G)\rightarrow\rep(G^\prime)$ corresponds to $(D(\varphi)^*)^\L$ under the equivalence.
    \end{enumerate}
\end{lem}

 In fact, $(\QCoh(G),\star)$ is equivalent to the monoidal category $(\QCoh(G)^\vee,(m^*)^\vee)$ naturally. Identifying $(\QCoh(G),m^*)$ comodules as $(\QCoh(G),\star)$ modules, we obtain:

\begin{cor}[Change of groups]\label{Ch_grp}
    There is an equivalence $\cC^{D(G)}\rightarrow\cC_G$, which is functorial on $G$ in the following sense. Suppose $G^\prime\rightarrow G$ is a morphism of groups and $D(G^\prime)$ acts on $\cC$, then $D(G)$ admits an induced action on $\cC$. Under the equivalences below, we have:
    \begin{center}\begin{tikzcd}
        \cC^{D(G)}& \cC_G \arrow[l, Rightarrow, no head] \arrow[d, "\oblv_*"']\\
        \cC^{D(G^\prime)} \arrow[u, "\oblv^*"'] & \cC_{G^\prime} \arrow[l, Rightarrow, no head] \arrow[lu, " ", phantom]
\end{tikzcd}\end{center}
    \begin{enumerate}
        \item If $G,G^\prime$ are discrete, then $\oblv^*\dashv\oblv_*$.
        \item If $G,G^\prime$ are diagonalizable, then $\oblv_*\dashv\oblv^*$.
    \end{enumerate}

\end{cor}

\begin{proof}
We prove (2) first. We calculate both sides as (co-)bar constructions
    $$\cC^{D(G)}=\lim\cb(\QCoh(D(G)),\cC);$$
    $$\cC_G=\colim\bc(\rep(G),\cC),$$
where the limit and colimit are taken over $\mathsf{Pr}^\L$. By \cite[Sec. D.3]{1aff}, taking right adjoints of the diagram $\bc(\rep(G),\cC)$ gives exactly rise to $\cb(\QCoh(D(G)),\cC)$. This completes the proof of (2).

Now we assume $G,G^\prime$ are discrete and prove (1). We define $\oblv^\star_D$ to be the composition $$\cC^{D(G)}=[(\cC^{D(G)})_{D(G^\prime)}]^{D(G^\prime)}\rightarrow[(\cC^{D(G)})_{D(G)}]^{D(G^\prime)}=\cC^{D(G^\prime)}$$ induced by the forgetful functor in (2) for $D(G)\rightarrow D(G^\prime)$. Similarly $\oblv_*^D$ is defined as the composition $$\cC_G=[(\cC_G)^{G^\prime}]_{G^\prime}\rightarrow[(\cC_G)^G]_{G^\prime}=\cC_{G^\prime}.$$
Then (2) induces an adjunction $\oblv_*^D\dashv\oblv^*_D$. To check $\oblv^*_D=\oblv_*$ and $\oblv_*^D=\oblv^*$, we express each equivariantization functor $(-)^G$ as $$(-)_{D(G)}=-\otimes_{\rep(D(G))}\Vect=-\otimes_{\QCoh(G)}\Vect$$ and apply Lem.\ref{dualgrp}.
\end{proof}

\begin{rmk}\label{eqn_calg}
    Recall that $\calg(\prlst)$ admits (small) limits, and the forgetful functor $\calg(\prlst)\rightarrow\prlst$ preserves limits. Therefore, if the co-Bar construction $\cb(\QCoh(G),\cC)$ is lifted into a cosimplicial object of $\calg(\prlst)$, then $\cC^G$ is symmetric monoidal and $\cC^G\rightarrow\cC$ is symmetric monoidal.
\end{rmk}

\subsection{Toric stacks}\label{torstk}

We treat a general case of toric stacks. If $X$ is a toric variety with torus $T$, then a morphism $G\rightarrow 
T$ between algebraic groups induces an action of $G$ on $X$. We say $[X/G]$ is a toric stack if $G$ is diagonalizable. To summarize, we define stacky fans and morphisms between them as follows.

\begin{defi}
    A {\bf stacky fan} is a tuple $(\sN,\Sigma, \sL,\beta)$, where \begin{enumerate}
        \item $\sN$ is a lattice;
        \item $\Sigma$ is a fan on $\sN_\RR$, as in \cite[Def. 3.1.2]{Cox};
        \item $\sL$ is a finitely generated abelian group;
        \item $\beta:\sM\rightarrow\sL$ is a morphism, where $\sM:=\sN^\vee$.
    \end{enumerate}
    For such a tuple, we assign a stack $\cX_{\Sigma,\beta}:=[X_\Sigma/D(\sL)]$, where $X_\Sigma$ is the toric variety associated to $(\sN,\Sigma)$ and $D(\sL)$ acts on $X_\Sigma$ by $D(\beta):D(\sL)\rightarrow D(\sM)$. Naturally a toric variety is regarded as a toric stack by $(\sN,\Sigma)=(\sN,\Sigma,0,0)$.
\end{defi}
\begin{defi}
    For stacky fans $(\sN,\Sigma, \sL,\beta)$ and $(\sN^\prime,\Sigma^\prime, \sL^\prime,\beta^\prime)$, a morphism consists of a fan morphism $\phi_\sN:(\sN,\Sigma)\rightarrow(\sN^\prime,\Sigma^\prime)$ and a morhism $\phi_\sL:\sL^\prime\rightarrow\sL$, satisfying $\phi_\sL\beta^\prime=\beta\phi_\sN^\vee$. Such a morphism induces a morphism $\phi$ between toric stacks as the composition
    $$[X_\Sigma/D(\sL)]\rightarrow[X_{\Sigma^\prime}/D(\sL)]\rightarrow[X_{\Sigma^\prime}/D(\sL^\prime)].$$
\end{defi}

\begin{rmk}\label{tor_stk_kuw}
In \cite{TorStk1,kuw16}, a stacky fan is encoded as $(L,\Sigma,N,\varphi)_{GS}$ where\begin{enumerate}
    \item Both $L,N$ are lattices;
    \item $\Sigma$ is a fan on $L_\RR$;
    \item $\varphi:L\rightarrow N$ is a morphism and $\coker\varphi$ is finite.
\end{enumerate}
In our notations, we have $$(L,\Sigma,N,\varphi)_{GS}=(L,\Sigma,(\ker\varphi)^\vee,\beta),$$ where $\beta=[\ker\varphi\rightarrow N]^\vee$. For a stacky fan $(\sN,\Sigma,\sL,\beta)$ with surjective $\beta$, we have$$(\sN,\Sigma,\sL,\beta)=(\sN,\Sigma,\coker(\beta^\vee),\varphi)_{GS},$$ where $\varphi:\sN\rightarrow\coker(\beta^\vee)$ is the natural map.
\end{rmk}

\subsection{Global quotients and covers}

We recall the $!$-tensor product $$-\otimes^!-:=\Delta^!(-\boxtimes-)$$ of ind-coherent sheaves is symmetric monoidal with the dualizing sheaf $\omega$ as unit. We always write it as $\otimes$ for simplicity. 

For an abelian Lie group $(T,e,m)$, there is a convolution product $$-\star-:=m_!(-\boxtimes-)$$ on $\Sh(T)$ with unit $k_e$.\footnote{When restricted to a skeleton $\Lambda$, the subcategory $\Sh_\Lambda(T)$ may fail to be closed with $\star$. Even if $\Sh_\Lambda(T)$ is symmetric monoidal with $\star$, the unit may differ from $k_e$.}

For a fan $\Sigma$ on $\sN_\RR$, we recall the FLTZ skeleton $$\Lambda_\Sigma:=\bigcup_{\sigma\in\Sigma}(\sigma^\perp/\sM)\times(-\sigma)\subseteq(\sM_\RR/\sM)\times\sN_\RR\cong T^*(\sM_\RR/\sM).$$ 

Kuwagaki proved the following theorem in \cite{kuw16}:

\begin{thm}\label{kuw} Here we restrict the results in \cite{kuw16} to toric varieties.
    \begin{enumerate}
        \item $\Sh_{\Lambda_\Sigma}(\sM_\RR/\sM)$ is symmetric monoidal under the convolution product;
        \item There is a symmetric monoidal equivalence $$K_\Sigma:(\indcoh(X_\Sigma),\otimes)\rightarrow(\Sh_{\Lambda_\Sigma}(\sM_\RR/\sM),\star);\footnote{This functor does not restrict to the functor $\kappa$ constructed in \cite{FLTZ,FLTZ09}, but they are related as in \cite[9.1(ii)]{kuw16}.}$$
        \item For an open inclusion $j$ induced by a subfan $\Sigma^\prime\subseteq\Sigma$, the equivalence interwines $j_*=(j^!)^\R$ and the inclusion $\Sh_{\Lambda_{\Sigma^\prime}}(\sM_\RR/\sM)\rightarrow\Sh_{\Lambda_\Sigma}(\sM_\RR/\sM)$;
        \item  For an smooth refinement $f:\Sigma^\prime\rightarrow\Sigma$, the equivalence interwines $f^!$ and the inclusion $\Sh_{\Lambda_\Sigma}(\sM_\RR/\sM)\rightarrow\Sh_{\Lambda_{\Sigma^\prime}}(\sM_\RR/\sM)$.
    \end{enumerate}
\end{thm}

Then we pass to quotient stacks by constructing the mirror as a (possibly not connected) cover of $\sM_\RR/\sM$. We need some preparation to ensure the (de-)equivariantization functors are compatible with the six-functor formalism of sheaves.

\begin{lem}[\'{E}tale descent of sheaves]\label{desc_shv}
    Suppose a finitely generated abelian group $G$ acts on a manifold $X$ freely and properly discontinuously ($X$ is a $G$-manifold for short), and $\pi:X\rightarrow X/G$ is the quotient. Let $\Lambda\subseteq T^*(X/G)$ be a closed conic subset and $\tilde{\Lambda}=\Lambda\times_{T^*(X/G)}T^*X\subseteq T^*X$. Let $G$ act on $\Sh_{\tilde{\Lambda}}(X)$ by translating the sheaves. Then $\pi^{-1}:\Sh_\Lambda(X/G)\rightarrow\Sh_{\tilde{\Lambda}}(X)$ factors through the forgetful functor $(\Sh_{\tilde{\Lambda}}(X))^G\rightarrow\Sh_{\tilde{\Lambda}}(X)$ and induces an equivalence $\Sh_\Lambda(X/G)\rightarrow(\Sh_{\tilde{\Lambda}}(X))^G$.
\end{lem}
For a proof, see for example \cite[Cor. 2.18]{CJ24}.


\begin{cor}[\'{E}tale descent of functors]\label{desc_of_fct}
    Suppose $\tilde{f}:X\rightarrow X^\prime$ is a $G$-equivariant map between $G$-manifolds, and $\Lambda\subseteq T^*(X/G);\Lambda^\prime\subseteq T^*(X^\prime/G)$. Then the inclusion $\tilde{f}^{-1}\Sh_{\tilde{\Lambda^\prime}}(X^\prime)\subseteq\Sh_{\tilde{\Lambda}}(X)$ holds iff $f^{-1}\Sh_{\Lambda^\prime}(X^\prime/G)\subseteq\Sh_\Lambda(X/G)$. In this case we have $f^{-1}=(\tilde{f}^{-1})^G$. Similar statements for $f^!,f_!,f_*$ hold.
\end{cor}

\begin{proof}

Without restrictions of singular supports, it is exactly \cite[Sec. 3.3]{eqvshvfun}.

To prove the general case, we assume the first inclusion ($\tilde{f}^{-1}\Sh_{\tilde{\Lambda^\prime}}(X^\prime)\subseteq\Sh_{\tilde{\Lambda}}(X)$, for instance) holds and take $(-)^G$ on the diagram below, where the vertical arrows are any pair of the functors. Dually taking $(-)_G$ proves another direction.

\begin{center}\begin{tikzcd}
\Sh_{\tilde{\Lambda}}(X) \arrow[r] \arrow[d, no head, dashed] & \Sh(X) \arrow[d, no head, dashed] \\
\Sh_{\tilde{\Lambda^\prime}}(X^\prime) \arrow[r]              & \Sh(X^\prime)                    
\end{tikzcd}\end{center}
\end{proof}


\begin{defi}\label{skeleton}
    For a stacky fan $(\sN,\Sigma,\sL,\beta)$, we define the stacky torus as the pushout $T_\beta=(\sM_\RR\times\sL)/\sM$ and identify $T_\beta/\sL=\sM_\RR/\sM$. Let $\Lambda_\beta\subseteq T^*T_\beta=T_\beta\times\sN_\RR$ be the preimage of $\Lambda_\Sigma\subseteq(\sM_\RR/\sM)\times\sN_\RR$, as in Lem.\ref{desc_shv}.
\end{defi}

\begin{thm}\label{main_thm_cat}
    There is a symmetric monoidal equivalence $$K_{\Sigma,\beta}:\indcoh(\cX_{\Sigma,\beta})\rightarrow\Sh_{\Lambda_\beta}(T_\beta).$$ In particular, it equals to $K_\Sigma$ in Thm.\ref{kuw} when $\sL=0$.
\end{thm}

\begin{proof}

We recall the constructions in \cite{FLTZ}. Set $$\Theta(\sigma)=\cO_{X_\sigma}\in\QCoh(X_\Sigma);$$
$$\Theta^\prime(\sigma)=p_!k_{\Int\sigma^\vee}[\rank\sN],$$
where $p:\sM_\RR\rightarrow\sM_\RR/\sM$.

If $\sigma\subseteq\sigma^\prime$, there is an identification $$\Hom(\Theta(\sigma^\prime),\Theta(\sigma))\cong k[\sM\cap\sigma^\vee]\cong\Hom(\Theta^\prime(\sigma^\prime),\Theta^\prime(\sigma)).$$ Under the $D(\sM)$-action on both sides, the following diagram commutes, where $\Delta$ is induced by the diagonal map $\sM\rightarrow\sM\times\sM$. Geometrically, a $k$-point $\varphi:\sM\rightarrow k^\times$ of $D(\sM)$ preserves the sheaves $\Theta,\Theta^\prime$ and acts on the morphisms by rescaling $\sM$.

\begin{center}\begin{tikzcd}
\Hom(\Theta(\sigma^\prime),\Theta(\sigma)) \arrow[r, Rightarrow, no head] \arrow[d] & k[\sM\cap\sigma^\vee] \arrow[r, Rightarrow, no head] \arrow[d, "\Delta"] & \Hom(\Theta^\prime(\sigma^\prime),\Theta^\prime(\sigma)) \arrow[d] \\
\Hom(\Theta(\sigma^\prime\times0),\Theta(\sigma\times0)) \arrow[r, Rightarrow, no head]           & k[\sM\cap\sigma^\vee][\sM] \arrow[r, Rightarrow, no head]                     & \Hom(\Theta^\prime(\sigma^\prime\times0),\Theta^\prime(\sigma\times0))          
\end{tikzcd}\end{center}

From the construction of $K_\Sigma$ in \cite{kuw16}, the discussion above shows $K_\Sigma$ is $D(\sM)$-equivariant. Therefore, $$K_{\Sigma,\beta}=(K_\Sigma)^{D(\sL)}:\indcoh(X_\Sigma)^{D(\sL)}\rightarrow\Sh_{\Lambda_\Sigma}(\sM_\RR/\sM)^{D(\sL)}=\Sh_{\Lambda_\Sigma}(\sM_\RR/\sM)_\sL$$ is the desired equivalence.

Finally, the category $\Sh_{\Lambda_\beta}(T_\beta)$ with the functor $K_{\Sigma,\beta}$ is enhanced to be symmetric monoidal by Lem.\ref{monoidal}.
\end{proof}

This result can be slightly generalized to a non-regular base.

\begin{cor}
    Suppose $k$ is a field and $S$ is a scheme locally of finite type over $k$. Then there is a symmetric monoidal equivalence 
    $$\indcoh(\cX_{\Sigma,\beta}\times_k S)\rightarrow\Sh_{\Lambda_\beta}(T_\beta,\indcoh(S)).$$
\end{cor}

\begin{proof}
    Apply Lurie tensor product on both sides.
\end{proof}

 In the rest of this section, we discuss the convolutional unit and closed structure of $\Sh_{\Lambda_\beta}(T_\beta)$, assuming Lem.\ref{monoidal} holds.

\begin{thm}\label{unit}
    Set $p:\sM_\RR\rightarrow T_\beta$. For any $\sigma\in\Sigma$, we define the dualizing sheaf $\1_\sigma:=k_{\Int\sigma^\vee}[\rank\sN]\in\Sh(\sM_\RR)$. When $\sigma$ varies, the sheaves form an inverse system, and the unit object is the limit $$\1_{\Sigma,\beta}:=p_!\lim_{\sigma\in\Sigma^{\op}}\1_\sigma.$$
\end{thm}

\begin{proof}
We assume $\Sigma$ is a single cone (with faces) first. The result in \cite{FLTZ} implies the case where $\beta$ is an isomorphism. Then apply Lem.\ref{fun_on_sch} to change the groups. Finally, the {\cech} resolution $\displaystyle\omega_{X_\Sigma}=\lim_{\sigma\in\Sigma^{\op}}\omega_{X_\sigma}$ completes the proof.
\end{proof}

\begin{rmk}
    For fixed $\beta:\sM\rightarrow\sL$, the sheaf $\1_{\Sigma,\beta}$ depends only on $|\Sigma|$. Indeed, for fans $\Sigma,\Sigma^\prime$ satisfying $|\Sigma|=|\Sigma^\prime|$, we select a common smooth refinement and apply Thm. \ref{kuw}(4).
\end{rmk}

\begin{prop}\label{closed}
    The functor $\hom^\star(x,y):=p_{1*}\hom(p_2^{-1}x,m^!y)$ restricts to a functor $\Sh_{\Lambda_\beta}(T_\beta)^{\op}\times\Sh_{\Lambda_\beta}(T_\beta)\rightarrow\Sh_{\Lambda_\beta}(T_\beta)$. Hence this is the closed structure of $\Sh_{\Lambda_\beta}(T_\beta)$.
\end{prop}

\begin{proof}
Because the functors are calculated \'{e}tale locally (\cite[Sec. 3.2]{eqvshvfun}), we may assume $\sL=0$. Then for $x,y\in\Sh_{\Lambda_\beta}(T_\beta)$, the decomposition $$\hom^\star(x,y)=\hom^\star(\1_\Sigma\star x,\1_\Sigma\star y)=\colim_{\sigma\in\Sigma}\lim_{\tau\in\Sigma^{op}}\hom^\star(\1_\sigma\star x,\1_{\tau}\star y)$$ implies that it suffices to prove $\ss[\hom^\star(\1_\sigma\star x,\1_\tau\star y)]\subseteq\Lambda_\Sigma$. By Thm.\ref{kuw}, $\ss[\1_\sigma\star x]\subseteq\Lambda_\sigma$ and $\ss[\1_\tau\star y]\subseteq\Lambda_\tau$. Finally \cite[Proof of Lem. 4.14]{kuw16} implies $\ss[\hom^\star(\1_\sigma\star x,\1_\tau\star y)]\subseteq\Lambda_{\sigma\cap\tau}$ and completes the proof.
\end{proof}

\section{Functoriality}
In this section we discuss how $K_{\Sigma,\beta}$ interacts with other functors.

\subsection{Left and right functoriality}
In this subsection, let $\phi:(\sN,\Sigma, \sL,\beta)\rightarrow(\sN^\prime,\Sigma^\prime, \sL^\prime,\beta^\prime)$ be a morphism between stacky fans.

\begin{defi}
    The map $\phi$ induces a morhism $\phi_T:T_{\beta^\prime}\rightarrow T_\beta$ between Lie groups. We say $\phi$ is {\bf left-functorial} if $\phi_{T!}\Sh_{\Lambda_{\beta^\prime}}(T_{\beta^\prime})\subseteq\Sh_{\Lambda_\beta}(T_\beta)$ and there is a natural equivalence $\phi_{T!}K_{\Sigma^\prime,\beta^\prime}\cong K_{\Sigma,\beta}\phi^!$ ($\phi^!\rightleftharpoons\phi_{T!}$ for short).
    Similarly, $\phi$ is {\bf right-functorial} if $(\phi^!)^\R\rightleftharpoons\phi_T^!$ holds.
\end{defi}

\begin{defi}
    For an inclusion of skeletons $\Lambda^\prime\subseteq\Lambda\subseteq T^*X$, the inclusion $$\iota_{\Lambda^\prime\rightarrow\Lambda}:\Sh_{\Lambda^\prime}(X)\rightarrow\Sh_\Lambda(X)$$ admits left and right adjoints, denoted by $\iota_{\Lambda^\prime\rightarrow\Lambda}^{\L,\R}$. They also satisfy \'{e}tale descent as in Cor.\ref{desc_of_fct}. We may omit the subscripts if the skeletons are clear from the context.
\end{defi}

\begin{lem}\label{fun_on_sch}
     If $\phi_{\mbox{scheme}}:(\sN,\Sigma)\rightarrow(\sN^\prime,\Sigma^\prime)$ is left-functorial, then $\phi$ is left-functorial. Similar statement holds for right-functoriality.
\end{lem}

\begin{proof}
    Factor $\phi$ as $(\sN,\Sigma,\sL,\beta)\xrightarrow{\phi_1}(\sN^\prime,\Sigma^\prime,\sL,\beta\circ\phi_\sM)\xrightarrow{\phi_2}(\sN^\prime,\Sigma^\prime,\sL^\prime,\beta^\prime)$. Cor.\ref{desc_of_fct} implies $\phi_1$ is left-functorial.
    
    For $\phi_2$, the assumptions of singular support hold automatically, hence we omit the subscripts here for simplicity. The exceptional pullback is a forgetful functor $\oblv^*:\indcoh(X_{\Sigma^\prime})^{D(\sL^\prime)}\rightarrow\indcoh(X_{\Sigma^\prime})^{D(\sL)}$, which is the left adjoint of the forgetful functor $\oblv_*:\Sh(\sM_\RR/\sM)_\sL\rightarrow\Sh(\sM_\RR/\sM)_{\sL^\prime}$ by Cor.\ref{Ch_grp}. And Lem.\ref{desc_shv} shows $\oblv_*=\pi^*:\Sh(T_{\beta\circ\phi_M})\rightarrow\Sh(T_{\beta^\prime})$. So $\phi_2^!\rightleftharpoons\pi_!$ as desired.
\end{proof}

\begin{lem}\label{cplt_to_sm}
    If $\sL=\sL^\prime=0$ and $(\sN,\Sigma)$ is smooth complete, then $\phi$ is left-functorial. 
\end{lem}
\begin{proof}
Factor $\phi_\sN$ into three maps $$c_\sN b_\sN a_\sN=\left(\begin{array}{cc}
    0 & \id_{\sN^\prime} 
\end{array}\right)\left(\begin{array}{cc}
    \id_\sN & 0 \\
    \phi_\sN & \id_{\sN^\prime}
\end{array}\right)\left(\begin{array}{c}
    \id_\sN \\
      0
\end{array}\right).$$ These are fan maps $(\sN,\Sigma)\xrightarrow{a}(\sN\oplus\sN^\prime,\Sigma\times0)\xrightarrow{b}(\sN\oplus\sN^\prime,\Sigma\times\Sigma^\prime)\xrightarrow{c}(\sN^\prime,\Sigma^\prime)$. 

We write $*=(0,0)$ for the trivial fan representing the toric variety $\spec k$. The map $*\rightarrow(\ZZ,0)$ is left-functorial by direct calculation. By K\"{u}nneth formula on both sides (\cite{kunneth}), $a$ is left-functorial. Because $\Sigma$ is complete, $k_e\in\Sh_{\Lambda_\Sigma}(\sM_\RR/\sM)$ is the unit. This proves $(\sN,\Sigma)\rightarrow*$ is left-functorial. Again by K\"{u}nneth formula, $c$ is left-functorial.

We refine $\Sigma\times\Sigma^\prime$ to $\Sigma_u$ such that $b_{\sN,\RR}(\Sigma\times0)$ is a subfan of $\Sigma_u$. Then we apply \cite[Thm. 11.1.9(b)]{Cox} to refine $\Sigma_u$ to a smooth fan $\Sigma_s$ such that $\Sigma_s$ contains all smooth cones in $\Sigma_u$. Hence we factor $b$ as $\Sigma\times0\xrightarrow{j}\Sigma_s\xrightarrow{p}\Sigma\times\Sigma^\prime$.

Recall that the inclusion of a subfan is right-functorial and a smooth refinement is left-functorial by Thm.\ref{kuw}. As $K_{\Sigma_s}$ is monoidal, we obtain $j^!\rightleftharpoons b_{T!}(-)\star(K\omega_{\Sigma\times0})$. Finally we conclude $ba$ is left-functorial, because $(ba)^!\rightleftharpoons a_{T!}(b_{T!}(-)\star(K\omega_{\Sigma\times0}))=a_{T!}b_{T!}(-)\star a_{T!}K\omega_{\Sigma\times0}=a_{T!}b_{T!}(-)\star Ka^!\omega_{\Sigma\times0}=a_{T!}b_{T!}$.
\end{proof}

\begin{thm}\label{adj_of_inc}
For any skeletons $\Lambda_\beta\subseteq\Pi_\beta\subseteq T^*T_\beta$ and $\Lambda_{\beta^\prime}\subseteq\Pi_{\beta^\prime}\subseteq T^*T_{\beta^\prime}$ such that $\phi_{T!}\Sh_{\Pi_{\beta^\prime}}(T_{\beta^\prime})\subseteq\Sh_{\Pi_\beta}(T_\beta)$, the relation $\phi^!\rightleftharpoons\iota_{\Lambda_\beta\rightarrow\Pi_\beta}^\L\phi_{T!}\iota_{\Lambda_{\beta^\prime}\rightarrow\Pi_{\beta^\prime}}$ holds. Equivalently, the relation $(\phi^!)^\R\rightleftharpoons\iota_{\Lambda_{\beta^\prime}\rightarrow\Pi_{\beta^\prime}}^\R\phi_T^!\iota_{\Lambda_\beta\rightarrow\Pi_\beta}$ holds under the condition $\phi_T^!\Sh_{\Pi_\beta}(T_\beta)\subseteq\Sh_{\Pi_{\beta^\prime}}(T_{\beta^\prime})$. The statement is independent when $\Pi_\beta,\Pi_{\beta^\prime}$ varies. In particular, there is a canonical choice $\Pi_\beta=T^*T_\beta;\Pi_{\beta^\prime}=T^*T_{\beta^\prime}$. As a consequence, $\phi$ is left-functorial iff $\phi_{T!}\Sh_{\Lambda_{\beta^\prime}}(T_{\beta^\prime})\subseteq\Sh_{\Lambda_\beta}(T_\beta)$ holds. Similar statement holds for right-functoriality.
\end{thm}

\begin{proof}
Compactify $\Sigma^\prime$ to a complete fan $\Sigma^\prime_c$. Then refine $\Sigma$ to a smooth fan $\Sigma_s$ and compactify $\Sigma_s$ to a smooth complete fan $\Sigma_c$ such that $\phi_{\sN}$ is again a fan morphism $(\sN,\Sigma_c)\rightarrow(\sN^\prime,\Sigma^\prime_c)$. We organize them into the following diagram.
\begin{center}\begin{tikzcd}
\Sigma \arrow[d, "\phi"]      & \Sigma_s \arrow[r, "b"] \arrow[l, "a"'] & \Sigma_c \arrow[d, "c"] \\
\Sigma^\prime \arrow[rr, "d"] && \Sigma^\prime_c        
\end{tikzcd}\end{center}

Now $b,d$ are right-functorial and $a,c$ are left functorial. We obtain $$\phi^!=(a^!)^\L a^!\phi^!d^!(d^!)^\R=(a^!)^\L b^!c^!(d^!)^\R\rightleftharpoons\iota_{\Sigma\rightarrow\Sigma_c}^\L\phi_{T!}\iota_{\Sigma^\prime\rightarrow\Sigma^\prime_c},$$ as desired.
\end{proof}

Here we recall the functorial estimation of singular supports to make use of this theorem.

\begin{lem}\cite[Prop. 5.4.4]{KS90}
    Suppose $f$ is a proper morphism between manifolds $f:Y\rightarrow X$. Differential of $f$ gives two maps $T^*Y\xleftarrow{f_d} T^*X\times_X Y\xrightarrow{f_\pi} T^*X$. For any $\cF\in\Sh(Y)$, the inclusion $\ss(f_*\cF)\subseteq f_\pi f_d^{-1}\ss(\cF)$ holds.
\end{lem}

\begin{lem}\label{monoidal}
    The map $\phi$ is left-functorial if $\phi_{\mbox{scheme}}$ is proper. In particular, picking $\phi_{\sN}=\id_\sN\oplus\id_\sN$ as the diagnoal morphism $\Delta:\cX_{\Sigma,\beta}\rightarrow\cX_{\Sigma,\beta}\times\cX_{\Sigma,\beta}$ proves that $\otimes\rightleftharpoons\star$, as in \cite[Cor 3.13]{FLTZ}.
\end{lem}

\begin{proof}
    Set $f=\phi_T$ as above and suppose $\sL=\sL^\prime=0$, then $$f_d=\id\times\phi_{\sN,\RR}:\sM^\prime_\RR/\sM^\prime\times\sN_\RR\rightarrow\sM^\prime_\RR/\sM^\prime\times\sN_\RR^\prime;$$
    $$f_\pi=\phi_T\times\id:\sM^\prime_\RR/\sM^\prime\times\sN_\RR\rightarrow\sM_\RR/\sM\times\sN_\RR.$$ For a cone $\sigma\in\Sigma^\prime$, we have $$f_\pi f_d^{-1}([\sigma^\perp/\sM^\prime]\times[-\sigma])=[\phi_{\sM,\RR}(\sigma^\perp)/\sM]\times[-\phi_{\sN,\RR}^{-1}\sigma].$$ If $\phi$ is proper, then this is a subset of $\Lambda_{\Sigma}$.
\end{proof}

\begin{exa}
    Let $\sN=\ZZ^2$ and $\Sigma$ is the fan representing $\AA^2$. The group $G=D(\ZZ/2)$ acts on $X_\Sigma$ encoded as $\sM=\ZZ^2\xrightarrow{(-1,1)}\ZZ/2=:\sL$. There is a canonical map $\phi:[X_\Sigma/G]\rightarrow X_\Sigma\git G$ from stacky quotient to GIT quotient. As the figures below, $\phi$ is induced by $\phi_\sN=\left(\begin{array}{cc}
    1 & 1 \\
    0 & 2
    \end{array}\right)$. We identify $T_\beta=\sM_\RR^\prime/\sM^\prime$ from the homeomorphism $\phi_T$, then $\phi^!\rightleftharpoons\iota$ because $\phi$ is left-functorial. This example is similar to \cite[Fig. 2]{FLTZ09}.
    
\noindent\begin{minipage}{0.5\linewidth}\begin{center}
    \begin{tikzpicture}
        \fill [color=lightgray] (0,0) -- (0,1.5) -- (1.5,1.5) -- (1.5,0);
        \draw [thick,<->] (0,1.5) -- (0,0) -- (1.5,0);
    \end{tikzpicture}
    \captionof{figure}{Fan of $\AA^2$}
\end{center}\end{minipage}
\begin{minipage}{0.5\linewidth}\begin{center}
    \begin{tikzpicture}
        \fill [color=lightgray] (0,0) -- (0.75,0) -- (0.75,1.5);
        \draw [thick,<->] (0.75,0) -- (0,0) -- (0.75,1.5);
        \node (B) at (0.75,1.5) [right]{$(1,2)$};
    \end{tikzpicture}
    \captionof{figure}{Fan of $\AA^2\git G$}
\end{center}\end{minipage}

\noindent\begin{minipage}{0.5\linewidth}\begin{center}
    \begin{tikzpicture}
        \draw (0,0) -- (2,0);
        \draw (0,2) -- (2,2);
        \draw [thick] (0,2) -- (2,1);
        \draw [thick] (0,1) -- (2,0);
        \draw [thick] (0,0) -- (0,2);
        \draw [thick] (2,0) -- (2,2);
        \foreach \y in {1,2,4,5}{
            \draw [thick,color=red] (0,0.333*\y) -- (-0.2,0.333*\y);
        }
        \foreach \y in {0,1}{
            \foreach \x in {1,2,3,4}{
                \draw [thick,color=red] (0.4*\x,\y+1-0.2*\x) -- (0.4*\x-0.1,\y+0.8-0.2*\x);
            }
        }
        \foreach \y in {1,2}{
            \fill [color=red!30] (0,\y) -- (-0.4,\y) -- (-0.15,\y-0.3);
            \draw [dashed,color=red] (-0.4,\y) -- (0,\y) -- (-0.15,\y-0.3);
        }
    \end{tikzpicture}
    \captionof{figure}{$T_\beta$ and $\Lambda_\beta$}
\end{center}\end{minipage}
\begin{minipage}{0.55\linewidth}\begin{center}
    \begin{tikzpicture}
        \draw (0,0) -- (2,0);
        \draw (0,2) -- (2,2);
        \draw [thick] (0,2) -- (2,1);
        \draw [thick] (0,1) -- (2,0);
        \draw [thick] (0,0) -- (0,2);
        \draw [thick] (2,0) -- (2,2);
        \foreach \y in {1,2,3,4,5}{
            \draw [thick,color=red] (0,0.333*\y) -- (-0.2,0.333*\y);
        }
        \foreach \y in {0,1}{
            \foreach \x in {0,1,2,3,4}{
                \draw [thick,color=red] (0.4*\x,\y+1-0.2*\x) -- (0.4*\x-0.1,\y+0.8-0.2*\x);
            }
        }
        \fill [color=red!30] (0,2) -- (-0.4,2) -- (-0.15,2-0.3);
        \draw [dashed,color=red] (-0.4,2) -- (0,2) -- (-0.15,2-0.3);
        
    \end{tikzpicture}
    \captionof{figure}{$\sM^\prime_\RR/\sM^\prime$ and {\color{red} $\Lambda_{\Sigma^\prime}$}}
\end{center}\end{minipage}

    More generally, for any proper morphism $\phi:(\sN,\Sigma)\rightarrow(\sN^\prime,\Sigma^\prime)$ such that $\phi_{\sN,\RR}$ is invertible, the discussion above applies to $$\tilde{\phi}:(\sN,\Sigma,\sL,\beta)\rightarrow(\sN^\prime,\Sigma^\prime),$$ where $$\beta:\sM\rightarrow\coker\phi_\sM=:\sL.$$ 
\end{exa}

\begin{thm}\label{lftr}
    We have $\phi^!\rightleftharpoons\1_{\Sigma,\beta}\star\phi_{T!}(-)$, where $\1_{\Sigma,\beta}$ is the convolutional unit of $\Sh_{\Lambda_\beta}(T_\beta)$. In particular, $\phi$ is left-functorial iff $\phi_{\mbox{scheme}}$ is proper.
\end{thm}

\begin{proof}
We compactify $\Sigma$ into a fan $\Sigma_p$ which is proper over $\Sigma^\prime$ and name the morphisms as $\Sigma\xrightarrow{j}\Sigma_p\xrightarrow{p}\Sigma^\prime$. In the smashing equivalence $$(j^!)^\R j^!=(-)\otimes(j^!)^\R j^!\omega_{\cX_{\Sigma_p,\beta}},$$ note that $j^!$ is monoidal and $(j^!)^\R\rightleftharpoons\iota_j$, this completes the proof of the first assertion. If $\phi$ is left-functorial, then $(j^!)^\R\phi^!=\phi_p^!$. Applying the equality to $\omega_{\cX_{\Sigma^\prime,\beta}}$, we find $(j^!)^\R$ preserves the unit, which happens only if $j=\id$.
\end{proof}

\begin{cor}\label{exc_funs}
The comparison $\phi_*\rightleftharpoons\iota^\L\phi_T^{-1}\iota$ holds if $\phi_\sL:\sL^\prime\rightarrow\sL$ is an isomorphism. By taking adjoints, we have $(\phi_*)^\R\rightleftharpoons\iota^\R\phi_{T*}\iota$ as well.
\end{cor}

\begin{proof}
    Compactify $\Sigma$ into a fan proper over $\Sigma^\prime$ as above. Then we apply Thm.\ref{adj_of_inc} to obtain $$\phi_*=p_*j_*=(p^!)^\L(j^!)^\R\rightleftharpoons (p_{T!}\iota)^\L\iota_j,$$ as desired.
\end{proof}

Similarly, the following condition implies right-functoriality.

\begin{defi}\cite{semistable}
    A morphism $\phi:(\sN,\Sigma)\rightarrow(\sN^\prime,\Sigma^\prime)$ between toric varieties is called {\bf weakly semistable} if it is flat and has reduced fibers. Equivalently, the following conditions hold for each cone $\sigma\in\Sigma$:\begin{enumerate}
        \item The image $\phi_{\sN,\RR}\sigma$ is a cone in $\Sigma^\prime$;
        \item The map $\phi_\sN$ restricts to a surjective map $\sigma\cap\sN\rightarrow(\phi_{\sN,\RR}\sigma)\cap\sN^\prime$.
    \end{enumerate}
\end{defi}

\begin{thm}\label{r_fct}
    If $\phi_{\mbox{scheme}}$ is weakly semistable, then $\phi$ is right-functorial. If furthermore $\phi_\sL$ is an isomorphism, then $\phi_*\rightleftharpoons\phi_T^{-1}$ holds as well.
\end{thm}

\begin{proof}
    By Lem.\ref{fun_on_sch}, we may assume $\sL=\sL^\prime=0$. For a sheaf $\cF\in\Sh_{\Lambda_\Sigma}(\sM_\RR/\sM)$, we have $\ss(\phi_T^!\cF)\subseteq\phi_T^\#\Lambda_\Sigma$ by \cite[Cor. 6.4.4]{KS90}. Following \cite[Prop. 6.2.4(iii)]{KS90}, we obtain 
    $$\phi_T^\#\Lambda_\Sigma=\bigcup_{\sigma\in\Sigma}[\phi_T^{-1}(\sigma^\perp/\sM)]\times[-\phi_{\sN,\RR}\sigma].$$
    Fix a cone $\sigma\in\Sigma$ and write $\tau:=\phi_{\sN,\RR}\sigma$. It suffices to prove $$\phi_T^{-1}(\sigma^\perp/\sM)\subseteq\tau^\perp/\sM^\prime,$$ or equivalently $$\phi_{\sM,\RR}^{-1}(\sigma^\perp+\sM)\subseteq\tau^\perp+\sM^\prime.$$
    By the assumption that $\phi$ is weakly semistable, we may choose a splitting $$\ZZ\sigma=\ZZ\tau\oplus\cK.$$
    Then we pick splittings of $\sN,\sN^\prime$ as follows.
    $$\sN=\ZZ\tau\oplus\cK\oplus\sN/\ZZ\sigma;$$
    $$\sN^\prime=\ZZ\tau\oplus\sN^\prime/\ZZ\tau.$$
    Under the dual splittings 
    $$\sM=(\ZZ\tau)^\vee\oplus\cK^\vee\oplus\ZZ\sigma^\perp;$$
     $$\sM^\prime=(\ZZ\tau)^\vee\oplus\ZZ\tau^\perp,$$
     the map $\phi_\sM$ is decomposed as $$\phi_{\sM}=\left(\begin{array}{cc}
         \id_{(\ZZ\tau)^\vee} & 0 \\
         0 & 0\\
         * & *
     \end{array}\right).$$ Finally, we conclude that $$\phi_{\sM,\RR}^{-1}(\sigma^\perp+\sM)=(\ZZ\tau)^\vee\oplus\tau^\perp=\tau^\perp+\sM^\prime.$$
\end{proof}

As an example, we investigate the toric stacks which are (represented by) toric varieties.

\begin{prop}\label{schematic}By \cite[Rem. 6.3]{TorStk2} and \cite[Cor. 6.5]{TorStk1}, we find that:
    \begin{enumerate}
        \item For a stacky fan $(\sN,\Sigma,\sL,\beta)$, the stack $\cX_{\Sigma,\beta}$ is a scheme iff $\sigma^\perp\cap\sM\xrightarrow{\beta}\sL$ is surjective for every (maximal) cone $\sigma\in\Sigma$;
        \item Suppose condition (1) holds, we write $\sK$ for the torsion-free quotient of $\coker[\sL^\vee\rightarrow\sN]$ and $\Phi:\sN\rightarrow\sK$. Then $\Phi_\RR(\sigma)$ is a cone in $\sK_\RR$ for every $\sigma\in\Sigma$. We organize them into a poset $\Sigma^\prime$. $\cX_{\Sigma,\beta}$ is a toric variety iff $\Phi_\RR:\Sigma\rightarrow\Sigma^\prime$ is an isomorphism between posets. In this case $\Sigma^\prime$ is a fan and the map $\Phi:(\sN,\Sigma,\sL,\beta)\rightarrow(\sK,\Sigma^\prime)$ induces an isomorphism between stacks;
        \item Suppose the condition in (2) holds. The map $\Phi:(\sN,\Sigma)\rightarrow(\sK,\Sigma^\prime)$ is weakly semistable.
    \end{enumerate}
\end{prop}
\begin{exa}
Suppose $\sN=\ZZ^2$ and $\Sigma$ is the fan representing $\AA^2$. The map $$\beta=\id\oplus(-\id):\sM=\ZZ^2\rightarrow\ZZ=\sL$$ defines a $\GG_m$-action on $X_\Sigma$ as the dashed line. The quotient $\cX=[X_\Sigma/\GG_m]$ is not represented by a scheme. We write $$\mathbf{n}=\left(\begin{array}{c}
     1\\
     -1 
\end{array}\right);$$
$$\phi_\sN:\sN\rightarrow\sN^\prime:=\sN/\mathbf{n};$$ $$\Sigma^\prime=\phi_{\sN,\RR}(\Sigma).$$ Then $\Sigma^\prime$ is a fan on $\sN^\prime_\RR$ representing $X_{\Sigma^\prime}=\AA^1$, and the map $$\phi_{\mbox{scheme}}:\AA^2\rightarrow\AA^1$$
$$(x,y)\mapsto xy$$
is weakly semistable. In the figure below, $T_\beta = \RR^2/\ZZ\mathbf{m}$ where $\mathbf{m}=(1,1)$ and the map $\phi_T:\sM^\prime_\RR/\sM^\prime\rightarrow T_\beta$ is an inclusion.

\noindent\begin{minipage}{0.5\linewidth}\begin{center}
    \begin{tikzpicture}
    \fill [color=lightgray] (0,0) -- (0,1.5) -- (1.5,1.5) -- (1.5,0);
    \draw [thick,<->] (0,1.5) -- (0,0) -- (1.5,0);
    \draw [dashed] (1,-1) -- (-1,1) node [above] {$G$};
\end{tikzpicture}
    \captionof{figure}{Fan of $\cX$}
\end{center}\end{minipage}
\begin{minipage}{0.5\linewidth}\begin{center}
    \begin{tikzpicture}
        \draw (0,2) -- (2,0);
        \draw (1,3) -- (3,1);
        \draw [dashed] (0,2) -- (1,3);
        \draw [dashed] (2,0) -- (3,1);
        \draw [thick] (1,3) -- (1,1) -- (3,1);
        \draw [thick] (0,2) -- (2,2) -- (2,0);
        \node (A) at (1,0.8) [below]{$0$};
        \foreach \x in {1,2}{
            \foreach \y in {1,2}{
            \fill[color=red!30] (\x,\y) -- (\x-0.3,\y) arc[start angle = 180, end angle = 270,radius=0.3] -- (\x,\y-0.3);
            }
        }
        \foreach \x in {0,1,...,4}{
            \foreach \y in {0,1}{
                \draw [thick,color=red] (0.4*\x+1-\y,\y+1) -- (0.4*\x+1-\y,\y+0.85);
                \draw [thick,color=red] (1+\y,0.4*\x+1-\y) -- (\y+0.85,0.4*\x+1-\y);
            }
        }
        \draw [thick,color=blue] (1,1) -- (2,2);
        \node (B) at (2,2) [right]{\color{blue} $\sM^\prime_\RR/\sM^\prime$};
    \end{tikzpicture}
\captionof{figure}{$T_\beta$ and {\color{red}$\Lambda_\beta$}}
\end{center}\end{minipage}
\end{exa}

\subsection{Application: Orbit closures}\label{closedorbits}
In this subsection we investigate the inclusion of an orbit closure. By equivariantization formalism, we restrict to non-stacky cases in this subsection. Similar comparison was established in \cite[Sec. 5]{toricdivisor} as Lagrangian correspondence.\\

In a fan $(\sN,\Sigma)$, we fix a cone $\tau\in\Sigma$ and the associated orbit closure $Y=\overline{O(\tau)}$.

\begin{defi}
    We set up some notations.\begin{enumerate}
        \item $\Star(\tau)=\{\sigma\in\Sigma|\tau\subseteq\sigma\}$;
        \item $\overline{\Star(\tau)}$ is the smallest subfan of $\Sigma$ containing $\Star(\tau)$. Let $U=X_{\overline{\Star(\tau)}}$;
        \item Let $p:\sN\rightarrow\sN/\ZZ\tau$ be the quotient map, then $p_\RR(\sigma)\subseteq\sN_\RR/\RR\tau$ is a cone for $\sigma\in\Star(\tau)$. We organize the image into a fan, and abuse notations to name it $\Star(\tau)$ again;
        \item We have $Y=X_{\Star(\tau)}$ and regard $p:U\rightarrow Y$ as a toric morphism.
        \end{enumerate}
\end{defi}

The mirror of $i_*$ along the closed immersion $i:Y\rightarrow U$ is computed as a convolution with a certain sheaf $\alpha_\tau$ preceded by a fully faithful pushforward.

\begin{thm}
    Set $\alpha_\tau=K_{\overline{\Star(\tau)}}i_*\omega_Y$. We have the following comparison. 
    $$i_*\rightleftharpoons \alpha_\tau\star p_{T!}(-);$$
    $$i^!\rightleftharpoons\iota_{\Lambda_{\Star(\tau)}}^\R p_T^!\hom^*(\alpha_\tau,-).$$
\end{thm}

\begin{proof}
    By projection formula and $p\circ i=\id$, we obtain
    $$i_*=(i^!)^\L=p^!(-)\otimes i_*\omega_Y;$$
    $$i^!=(p^!)^\R\hom(i_*\omega_Y,-).$$
    Then apply Prop.\ref{closed}, Thm.\ref{adj_of_inc}, Lem.\ref{monoidal} and Thm.\ref{lftr}.
\end{proof}

\begin{rmk}
    If $\overline{\Star(\tau)}$ is simplicial, then $p$ is weakly semistable so we have $i^!\rightleftharpoons p_T^!\hom^*(\alpha_\tau,-)$. This comparison was formulated as microlocalization in \cite{shende-gammage} for smooth cases.
\end{rmk}

\begin{exa}
    Let $\Sigma$ be the fan representing the canonical bundle of $\PP^1$ and $Y$ is the zero section. Up to an invertible sheaf, the sheaf $i_*\omega_Y$ is calculated as $$\cF:=\cofib[\omega_X(-Y)\rightarrow\omega_X].$$
    We lift the sheaves to the stack $[X/\GG_m^2]$ and apply the calculation in \cite{FLTZ}. In the plane $\RR_{x,y}=\sM_\RR$, we set three open half-planes
    $$H_1=\{x>1\};H_2=\{x>y\};H_3=\{x>-y\},$$
    and write $H_{ij}=H_i\cap H_j$. From the {\cech} resolution of $\omega_X(-Y)$, we have $$K\omega_X(-Y)=\fib[k_{H_{12}}\oplus k_{H_{13}}\rightarrow k_{H_1}][2]=k_{H_{123}}[2].$$ Then $$K\cF=\cofib[k_{H_{123}}\rightarrow k_{H_{23}}][2]=k_{P}[2],$$
    where $P=H_{23}-H_{123}$.\\
    \noindent\begin{minipage}{0.5\linewidth}\begin{center}
    \begin{tikzpicture}
    \fill[color=lightgray](0,0) -- (1.5,-1.5) -- (1.5,1.5);
    \foreach \y in {-1.5,0,1.5}{
    \draw[thick,->] (0,0) -- (1.5,\y);
    }
    \node (A) at (1.5,0) [right] {$\tau$};
\end{tikzpicture}
    \captionof{figure}{The fan $\Sigma$}
\end{center}\end{minipage}
\begin{minipage}{0.5\linewidth}\begin{center}
    \begin{tikzpicture}
        \fill [color=lightgray] (0,0) -- (1,-1) --(1,1);
        \draw [dashed] (0,0) -- (1,1);
        \draw [thick] (1,1) -- (2,2);
        \draw [dashed] (0,0) -- (1,-1);
        \draw [thick] (1,-1) -- (2,-2);
        \draw [thick] (1,2) -- (1,-2);
        \node (B) at (1,2) [right]{$H_1$};
        \node (C) at (2,2) [below]{$H_2$};
        \node (D) at (2,-2) [above]{$H_3$};
        \draw [<->] (2,0) -- (0,0) -- (0,2);
        \node (E) at (2,0) [above]{$x$};
        \node (F) at (0,2) [right]{$y$};
    \end{tikzpicture}
    \captionof{figure}{$P\subseteq\sM_\RR$}
\end{center}\end{minipage}

\end{exa}

\printbibliography

@book {HTT,
    AUTHOR = {Lurie, Jacob},
     TITLE = {Higher topos theory},
    SERIES = {Annals of Mathematics Studies},
    VOLUME = {170},
 PUBLISHER = {Princeton University Press, Princeton, NJ},
      YEAR = {2009},
     PAGES = {xviii+925},
      ISBN = {978-0-691-14049-0; 0-691-14049-9},
   MRCLASS = {18-02 (18B25 18E35 18G30 18G55 55U40)},
  MRNUMBER = {2522659},
MRREVIEWER = {Mark\ Hovey},
       DOI = {10.1515/9781400830558},
       URL = {https://doi.org/10.1515/9781400830558},
}

@book {HA,
  title={Higher algebra},
  author={Lurie, Jacob},
  year={2017},
  publisher={Harvard University}
}

@article {lemma,
    AUTHOR = {Robalo, Marco and Schapira, Pierre},
     TITLE = {A lemma for microlocal sheaf theory in the
              {$\infty$}-categorical setting},
   JOURNAL = {Publ. Res. Inst. Math. Sci.},
  FJOURNAL = {Publications of the Research Institute for Mathematical
              Sciences},
    VOLUME = {54},
      YEAR = {2018},
    NUMBER = {2},
     PAGES = {379--391},
      ISSN = {0034-5318,1663-4926},
   MRCLASS = {55U35 (18F99 32C38 35A27 55P42)},
  MRNUMBER = {3784874},
MRREVIEWER = {Birgit\ Richter},
       DOI = {10.4171/PRIMS/54-2-5},
       URL = {https://doi.org/10.4171/PRIMS/54-2-5},
}

@article {unbounded,
    AUTHOR = {Spaltenstein, N.},
     TITLE = {Resolutions of unbounded complexes},
   JOURNAL = {Compositio Math.},
  FJOURNAL = {Compositio Mathematica},
    VOLUME = {65},
      YEAR = {1988},
    NUMBER = {2},
     PAGES = {121--154},
      ISSN = {0010-437X,1570-5846},
   MRCLASS = {18E25 (32C35)},
  MRNUMBER = {932640},
MRREVIEWER = {Michael\ M.\ Kapranov},
       URL = {http://www.numdam.org/item?id=CM_1988__65_2_121_0},
}

@book {KS90,
    AUTHOR = {Kashiwara, Masaki and Schapira, Pierre},
     TITLE = {Sheaves on manifolds},
    SERIES = {Grundlehren der mathematischen Wissenschaften [Fundamental
              Principles of Mathematical Sciences]},
    VOLUME = {292},
      NOTE = {With a chapter in French by Christian Houzel},
 PUBLISHER = {Springer-Verlag, Berlin},
      YEAR = {1990},
     PAGES = {x+512},
      ISBN = {3-540-51861-4},
   MRCLASS = {58G07 (18F20 32C38 35A27)},
  MRNUMBER = {1074006},
MRREVIEWER = {Michael\ M.\ Kapranov},
       DOI = {10.1007/978-3-662-02661-8},
       URL = {https://doi.org/10.1007/978-3-662-02661-8},
}

@incollection {1aff,
    AUTHOR = {Gaitsgory, Dennis},
     TITLE = {Sheaves of categories and the notion of 1-affineness},
 BOOKTITLE = {Stacks and categories in geometry, topology, and algebra},
    SERIES = {Contemp. Math.},
    VOLUME = {643},
     PAGES = {127--225},
 PUBLISHER = {Amer. Math. Soc., Providence, RI},
      YEAR = {2015},
      ISBN = {978-1-4704-1557-0},
   MRCLASS = {14F05 (14A20)},
  MRNUMBER = {3381473},
MRREVIEWER = {Jon\ Eivind\ Vatne},
       DOI = {10.1090/conm/643/12899},
       URL = {https://doi.org/10.1090/conm/643/12899},
}

@article {kuw16,
    AUTHOR = {Kuwagaki, Tatsuki},
     TITLE = {The nonequivariant coherent-constructible correspondence for
              toric stacks},
   JOURNAL = {Duke Math. J.},
  FJOURNAL = {Duke Mathematical Journal},
    VOLUME = {169},
      YEAR = {2020},
    NUMBER = {11},
     PAGES = {2125--2197},
      ISSN = {0012-7094,1547-7398},
   MRCLASS = {53D37 (14A20 14J33 14M25 35A27)},
  MRNUMBER = {4132582},
MRREVIEWER = {Amin\ Gholampour},
       DOI = {10.1215/00127094-2020-0011},
       URL = {https://doi.org/10.1215/00127094-2020-0011},
}

@book {eqvshvfun,
    AUTHOR = {Bernstein, Joseph and Lunts, Valery},
     TITLE = {Equivariant sheaves and functors},
    SERIES = {Lecture Notes in Mathematics},
    VOLUME = {1578},
 PUBLISHER = {Springer-Verlag, Berlin},
      YEAR = {1994},
     PAGES = {iv+139},
      ISBN = {3-540-58071-9},
   MRCLASS = {55N91 (14M25 18E30 54B40 55N30)},
  MRNUMBER = {1299527},
MRREVIEWER = {Yi\ Hu},
       DOI = {10.1007/BFb0073549},
       URL = {https://doi.org/10.1007/BFb0073549},
}

@misc{bai,
      title={Toric Mirror Symmetry for Homotopy Theorists}, 
      author={Qingyuan Bai and Yuxuan Hu},
      year={2025},
      eprint={2501.06649},
      archivePrefix={arXiv},
      primaryClass={math.AG},
      url={https://arxiv.org/abs/2501.06649}, 
}

@misc{kunneth,
      title={Duality, K\"unneth formulae, and integral transforms in microlocal geometry}, 
      author={Christopher Kuo and Wenyuan Li},
      year={2024},
      eprint={2405.15211},
      archivePrefix={arXiv},
      primaryClass={math.SG},
      url={https://arxiv.org/abs/2405.15211}, 
}

@book {Cox,
    AUTHOR = {Cox, David A. and Little, John B. and Schenck, Henry K.},
     TITLE = {Toric varieties},
    SERIES = {Graduate Studies in Mathematics},
    VOLUME = {124},
 PUBLISHER = {American Mathematical Society, Providence, RI},
      YEAR = {2011},
     PAGES = {xxiv+841},
      ISBN = {978-0-8218-4819-7},
   MRCLASS = {14M25 (05A15 05E45 52B12)},
  MRNUMBER = {2810322},
MRREVIEWER = {Ivan\ Arzhantsev},
       DOI = {10.1090/gsm/124},
       URL = {https://doi.org/10.1090/gsm/124},
}

@article {TorStk1,
    AUTHOR = {Geraschenko, Anton and Satriano, Matthew},
     TITLE = {Toric stacks {I}: {T}he theory of stacky fans},
   JOURNAL = {Trans. Amer. Math. Soc.},
  FJOURNAL = {Transactions of the American Mathematical Society},
    VOLUME = {367},
      YEAR = {2015},
    NUMBER = {2},
     PAGES = {1033--1071},
      ISSN = {0002-9947,1088-6850},
   MRCLASS = {14D23 (14M25)},
  MRNUMBER = {3280036},
MRREVIEWER = {Fabio\ Perroni},
       DOI = {10.1090/S0002-9947-2014-06063-7},
       URL = {https://doi.org/10.1090/S0002-9947-2014-06063-7},
}

@article {TorStk2,
    AUTHOR = {Geraschenko, Anton and Satriano, Matthew},
     TITLE = {Toric stacks {II}: {I}ntrinsic characterization of toric
              stacks},
   JOURNAL = {Trans. Amer. Math. Soc.},
  FJOURNAL = {Transactions of the American Mathematical Society},
    VOLUME = {367},
      YEAR = {2015},
    NUMBER = {2},
     PAGES = {1073--1094},
      ISSN = {0002-9947,1088-6850},
   MRCLASS = {14D23 (14M25)},
  MRNUMBER = {3280037},
MRREVIEWER = {Fabio\ Perroni},
       DOI = {10.1090/S0002-9947-2014-06064-9},
       URL = {https://doi.org/10.1090/S0002-9947-2014-06064-9},
}

@article {shende,
    AUTHOR = {Shende, Vivek},
     TITLE = {Toric mirror symmetry revisited},
   JOURNAL = {C. R. Math. Acad. Sci. Paris},
  FJOURNAL = {Comptes Rendus Math\'ematique. Acad\'emie des Sciences. Paris},
    VOLUME = {360},
      YEAR = {2022},
     PAGES = {751--759},
      ISSN = {1631-073X,1778-3569},
   MRCLASS = {14J33 (14A20 14M25)},
  MRNUMBER = {4449875},
MRREVIEWER = {Lei\ Yang},
       DOI = {10.5802/crmath.304},
       URL = {https://doi.org/10.5802/crmath.304},
}

@book {vol1,
    AUTHOR = {Gaitsgory, Dennis and Rozenblyum, Nick},
     TITLE = {A study in derived algebraic geometry. {V}ol. {I}.
              {C}orrespondences and duality},
    SERIES = {Mathematical Surveys and Monographs},
    VOLUME = {221},
 PUBLISHER = {American Mathematical Society, Providence, RI},
      YEAR = {2017},
     PAGES = {xl+533},
      ISBN = {978-1-4704-3569-1},
   MRCLASS = {14F05 (18D05 18G55)},
  MRNUMBER = {3701352},
MRREVIEWER = {Adrian\ Langer},
       DOI = {10.1090/surv/221.1},
       URL = {https://doi.org/10.1090/surv/221.1},
}

@article {FLTZ,
    AUTHOR = {Fang, Bohan and Liu, Chiu-Chu Melissa and Treumann, David and
              Zaslow, Eric},
     TITLE = {A categorification of {M}orelli's theorem},
   JOURNAL = {Invent. Math.},
  FJOURNAL = {Inventiones Mathematicae},
    VOLUME = {186},
      YEAR = {2011},
    NUMBER = {1},
     PAGES = {79--114},
      ISSN = {0020-9910,1432-1297},
   MRCLASS = {14F05 (14M25)},
  MRNUMBER = {2836052},
MRREVIEWER = {Pawel\ Sosna},
       DOI = {10.1007/s00222-011-0315-x},
       URL = {https://doi.org/10.1007/s00222-011-0315-x},
}

@article {shende-gammage,
    AUTHOR = {Gammage, Benjamin and Shende, Vivek},
     TITLE = {Mirror symmetry for very affine hypersurfaces},
   JOURNAL = {Acta Math.},
  FJOURNAL = {Acta Mathematica},
    VOLUME = {229},
      YEAR = {2022},
    NUMBER = {2},
     PAGES = {287--346},
      ISSN = {0001-5962,1871-2509},
   MRCLASS = {14J33 (14F08)},
  MRNUMBER = {4554224},
MRREVIEWER = {Guangbo\ Xu},
       DOI = {10.4310/acta.2022.v229.n2.a2},
       URL = {https://doi.org/10.4310/acta.2022.v229.n2.a2},
}

@article {FLTZ09,
    AUTHOR = {Fang, Bohan and Liu, Chiu-Chu Melissa and Treumann, David and
              Zaslow, Eric},
     TITLE = {The coherent-constructible correspondence for toric
              {D}eligne-{M}umford stacks},
   JOURNAL = {Int. Math. Res. Not. IMRN},
  FJOURNAL = {International Mathematics Research Notices. IMRN},
      YEAR = {2014},
    NUMBER = {4},
     PAGES = {914--954},
      ISSN = {1073-7928,1687-0247},
   MRCLASS = {14D23 (14F05)},
  MRNUMBER = {3168399},
MRREVIEWER = {Hsian-Hua\ Tseng},
       DOI = {10.1093/imrn/rns235},
       URL = {https://doi.org/10.1093/imrn/rns235},
}

@article {WrpShv,
    AUTHOR = {Kuo, Christopher},
     TITLE = {Wrapped sheaves},
   JOURNAL = {Adv. Math.},
  FJOURNAL = {Advances in Mathematics},
    VOLUME = {415},
      YEAR = {2023},
     PAGES = {Paper No. 108882, 71},
      ISSN = {0001-8708,1090-2082},
   MRCLASS = {18F20 (32C38 53D37 53D40 55N30)},
  MRNUMBER = {4543073},
MRREVIEWER = {Mee\ Seong\ Im},
       DOI = {10.1016/j.aim.2023.108882},
       URL = {https://doi.org/10.1016/j.aim.2023.108882},
}

@inproceedings{Bondal,
  title={Derived categories of toric varieties},
  author={Bondal, Alexey},
  booktitle={Convex and Algebraic geometry, Oberwolfach conference reports, EMS Publishing House},
  volume={3},
  pages={284--286},
  year={2006}
}

@misc{Tr10,
      title={Remarks on the nonequivariant coherent-constructible correspondence for toric varieties}, 
      author={David Treumann},
      year={2010},
      eprint={1006.5756},
      archivePrefix={arXiv},
      primaryClass={math.AG},
      url={https://arxiv.org/abs/1006.5756}, 
}

@article {SS,
    AUTHOR = {Scherotzke, Sarah and Sibilla, Nicol\`o},
     TITLE = {The non-equivariant coherent-constructible correspondence and
              a conjecture of {K}ing},
   JOURNAL = {Selecta Math. (N.S.)},
  FJOURNAL = {Selecta Mathematica. New Series},
    VOLUME = {22},
      YEAR = {2016},
    NUMBER = {1},
     PAGES = {389--416},
      ISSN = {1022-1824,1420-9020},
   MRCLASS = {14M25 (14F05 32S60 53D37)},
  MRNUMBER = {3437841},
MRREVIEWER = {Colin\ Diemer},
       DOI = {10.1007/s00029-015-0193-y},
       URL = {https://doi.org/10.1007/s00029-015-0193-y},
}

@article {CJ24,
    AUTHOR = {Clausen, Dustin and Jansen, Mikala \O rsnes},
     TITLE = {The reductive {B}orel-{S}erre compactification as a model for
              unstable algebraic {K}-theory},
   JOURNAL = {Selecta Math. (N.S.)},
  FJOURNAL = {Selecta Mathematica. New Series},
    VOLUME = {30},
      YEAR = {2024},
    NUMBER = {1},
     PAGES = {Paper No. 10, 93},
      ISSN = {1022-1824,1420-9020},
   MRCLASS = {18F25 (18N60 19D06 32S60)},
  MRNUMBER = {4683160},
MRREVIEWER = {Daniel\ A.\ Ramras},
       DOI = {10.1007/s00029-023-00900-8},
       URL = {https://doi.org/10.1007/s00029-023-00900-8},
}

@article {semistable,
    AUTHOR = {Molcho, Sam},
     TITLE = {Universal stacky semistable reduction},
   JOURNAL = {Israel J. Math.},
  FJOURNAL = {Israel Journal of Mathematics},
    VOLUME = {242},
      YEAR = {2021},
    NUMBER = {1},
     PAGES = {55--82},
      ISSN = {0021-2172,1565-8511},
   MRCLASS = {14A20 (14A21 14M25)},
  MRNUMBER = {4282076},
MRREVIEWER = {Damiano\ Fulghesu},
       DOI = {10.1007/s11856-021-2118-0},
       URL = {https://doi.org/10.1007/s11856-021-2118-0},
}

@article {toricdivisor,
    AUTHOR = {Hanlon, A. and Hicks, J.},
     TITLE = {Aspects of functoriality in homological mirror symmetry for
              toric varieties},
   JOURNAL = {Adv. Math.},
  FJOURNAL = {Advances in Mathematics},
    VOLUME = {401},
      YEAR = {2022},
     PAGES = {Paper No. 108317, 92},
      ISSN = {0001-8708,1090-2082},
   MRCLASS = {53D37 (14J33 14M25 14T20)},
  MRNUMBER = {4394684},
MRREVIEWER = {Helge\ Ruddat},
       DOI = {10.1016/j.aim.2022.108317},
       URL = {https://doi.org/10.1016/j.aim.2022.108317},
}
\end{document}